\documentclass{article}
\author{J\"org Peters, University of Florida}
\title{On $G^1$ stitched bi-cubic B\'ezier patches with arbitrary topology}
\usepackage{times}
\usepackage{amsmath}
\usepackage{amssymb}
\usepackage{amsfonts}
\usepackage{epsfig}
\usepackage[usenames,dvipsnames]{color}
\usepackage{subfigure}
\usepackage{overpic}

\newcommand{\bp}{\mathbf{p}}
\newcommand{\bm}{\mathbf{m}}
\newcommand{\bq}{\mathbf{q}}
\newcommand{\bv}{\mathbf{v}}
\newcommand{\bw}{\mathbf{w}}

\newcommand{\msh}{{\cal{M}}}
\newcommand{\R}{{\mathbb{R}}}

\newcommand{\figref}[1]{Fig.~\ref{#1}}
\newcommand{\secref}[1]{Section~\ref{#1}}

\newcommand{\thmref}[1]{Theorem~\ref{#1}}
\newtheorem{theorem}{Theorem}

\begin{document}
\maketitle
\begin{abstract}
Lower bounds on the generation of smooth bi-cubic surfaces
imply that geometrically smooth ($G^1$) constructions
need to satisfy conditions on the connectivity and layout.
In particular, quadrilateral meshes of arbitrary topology can not in general 
be covered with $G^1$-connected B\'ezier patches of bi-degree 3
using the layout proposed in \cite{Akle:2017:IM}.
This paper analyzes whether the pre-refinement of the input mesh by repeated
Doo-Sabin subdivision proposed in that paper yields an exception.

\end{abstract}
\section{Introduction}


For many applications, for example artistic rendering and 
sculpting, a few subdivision steps provide a pleasing rounding
of the original polyhedral shape. 
The simplicity of subdivision with small, local stencils (refinement rules)
is appealing and, in particular Catmull-Clark subdivision
\cite{Catmull-1978-CC} is a staple of geometric modeling 
environments for creating computer graphics assets.
However Catmull-Clark subdivision has also been demonstrated to lead to 
shape deficiencies, such as pinching of highlight lines,
that can be traced back to the simple stencil-based rules 
\cite{Karciauskas:2004:SCS,Karciauskas:ISR:2016}.

The algorithm of \cite{Akle:2017:IM} proposes an approach to obtaining 
`$C^2$ continuous Bi-Cubic B\'ezier patches that are guaranteed to be stitched
with $G^1$ continuity regardless of the underlying mesh topology'.
This approach consists of applying  not Catmull-Clark but
\emph{Doo-Sabin subdivision} to an initial polyhedral input mesh.
The approach then derives quadrilateral facets and B\'ezier control 
points from the refined mesh and constructs $n$ bi-cubic patches for each 
$n$-sided facet.

Beyond demonstrating pleasant rounding,
\cite{Akle:2017:IM} emphasizes that the result is a `smooth surface with 
$G^1$ continuity'
\footnote{$G^1$ is typeset as $G_1$ in several places in \cite{Akle:2017:IM}.
}.
%
If true, this result would be remarkable.
It would contradict the restrictions on bi-cubic $G^1$ 
spline complexes that \cite[Section 3]{Peters:2009:CSS} 
derived and that prompted the special constructions 
in \cite{Peters:2015:PSS,Sarov:RPG:2016}.
If \cite{Akle:2017:IM} were correct then these special constructions
published earlier in the same conference series would be superfluous!

Below we show that, while the surfaces generated by the approach of
\cite{Akle:2017:IM} often appear to be smooth, in general they are not.

\noindent{\bf Overview.}
\secref{sec:thm} summarizes the algorithm in \cite{Akle:2017:IM}
and the lower bound result of \cite{Peters:2009:CSS}as it pertains to 
bi-cubic $G^1$ constructions.
\secref{sec:counter} provides an explicit, minimal counterexample 
to the claim that the approach in \cite{Akle:2017:IM} generates
$G^1$ surfaces.
\secref{sec:alternative} discusses options for constructing 
both formally smooth and near-smooth bi-3 constructions.


\section{$G^1$ continuity, 
the construction of \cite{Akle:2017:IM} and a theorem}
\label{sec:thm}
\newcommand{\vloc}{vertex-localized}
\newcommand{\ik}{edge knot}
The construction of \cite{Akle:2017:IM}
applies two steps of Doo-Sabin subdivision
to an initial polyhedral input mesh $\msh$ and then places 
the corners of bicubic patches at the Doo-Sabin limit points of
the facets obtained in the initial subdivision 
(Fig 5 of \cite{Akle:2017:IM}).
That is every vertex and every face of $\msh$ has a corner 
of a bi-3 patch associated with it.
This layout looks more general, and therefore more challenging
than the one in \cite{conf/gmp/HahmannBC08} which assumed 
that the input mesh has quadrilateral faces
and used $2\times 2$ bi-cubics to cover them.

Denote by $\bv$ and $\bw$ limit points associated with  
adjacent facets of $\msh$ (see \figref{fig:tet}).
Since $\msh$ is unrestricted and 
$\bv$ and $\bw$ and their tangent planes can be freely adjusted.
This independence is typically desirable for flexibility of modeling.
The construction in \cite{Akle:2017:IM}
therefore $G^1$ \emph{\vloc} in the sense that the
Taylor expansion at $\bv$ is not tightly linked to that at $\bw$.
Since it does not matter in the construction
whether $\bv$ or $\bw$ is listed first,
the construction along the common boundary
is also \emph{unbiased}, and this is also typically desirable.
The \emph{unbiased $G^1$ constraints} between two patches
$\bp, \bq: (u,v) \to \R^3$ along $\bp(u,0) = \bq(u,0)$ are
\begin{equation}
   \partial_2\bp(u,0) +\partial_2\bq(u,0)
   =
   \alpha(u) \partial_1 \bp(u,0). \label{eq:g1}
\end{equation}
When, at the split point $\bm$ of 
the edge $\bv,\bw$, the four bicubic patches join $C^1$ 
(see e.g.\ Fig 8 of \cite{Akle:2017:IM})
the following theorem applies.
\begin{theorem} [\cite{Peters:2009:CSS}: two double \ik s needed]
In general, using splines of degree bi-3 for
a \vloc\ unbiased $G^1$ construction
without forced linear boundary segments,
the splines must have at least two internal double knots.
\label{thm:twodblknt}
\end{theorem}
In other words, \thmref{thm:twodblknt} states that to satisfy $G^1$ constraints
along $\bv$ and $\bw$
(and not have straight line segments embedded in the surface), 
three rather than the constructed two polynomial 
boundary segments are needed to connect $\bv$ and $\bw$.
One might hope that the initialization via Doo-Sabin or 
adding or leaving out some adjustment does side-step the 
assumptions of \thmref{thm:twodblknt}.
The next section therefore
looks more closely at the construction of \cite{Akle:2017:IM}.

Below the bicubic tensor-product polynomial surface patches $\bp$, $\bq$ 
of bi-degree $3$ are in the following expressed in Bernstein-B\'ezier (BB) form,
e.g.\
\[
   \bp(u,v):=\sum_{i=0}^3\sum_{j=0}^3\bp_{ij}B^3_i(u)B^3_j(v),
      \quad (u,v) \in \square:=[0..1]^2,
\]
where $B^3_k(t) := \binom{3}{k}(1-t)^{3-k}t^k$ is
the Bernstein-B\'ezier (BB) polynomials of degree $3$
and $\bp_{ij} \in \R^3$ are the BB coefficients \cite{Farin02,Prautzsch02}.


\section{A Counterexample: an input mesh where \cite{Akle:2017:IM} 
does not yield a $G^1$ output}
\label{sec:counter}
\def\wid{0.6\linewidth}
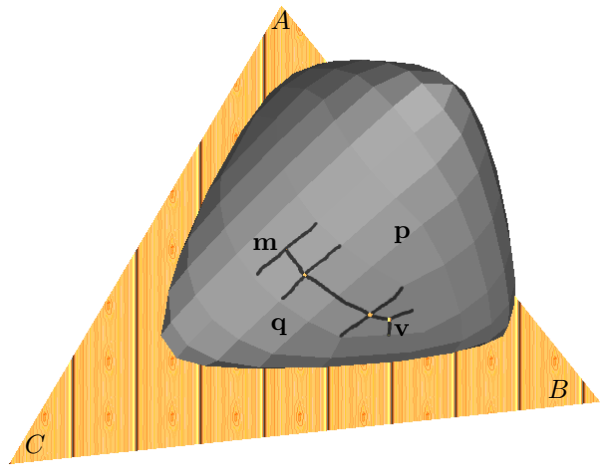
\begin{figure}[h]
   \centering
   \begin{overpic}[scale=.5,tics=10]{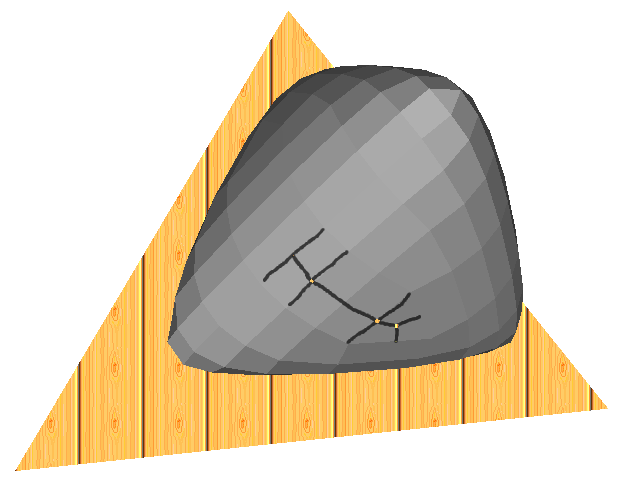} 
   \put(42,38){$\bm$}
   \put(65,40){$\bp$}
   \put(45,25){$\bq$}
   \put(65,24){$\bv$}
   \put(90,14){$B$}
   \put(5,5){$C$}
   \put(45,74){$A$}
   \end{overpic}
\label{fig:tet}
\caption{
Counterexample: the input is regular tetrahedron, only one of whose 
faces is shown with a wood texture. 
The grey quad-mesh is the result of applying three steps of
Doo-Sabin subdivision. The subnet of 12 Bernstein-B\'ezier
control points of interest are sketched on the 
refined mesh: from the 4-valent point $\bm$ to the 3-valent point
$\bv$, these are the BB-coefficients of \eqref{eq:g1coef} that influence the 
$G^1$ continuity between the two bi-3 patches $\bp$ and $\bq$.
}
\end{figure}

Since the algorithm of \cite{Akle:2017:IM} applies initially multiple
steps of Doo-Sabin subdivision, the challenge of finding a simple
explicit counterexample seems formidable.
Yet, the simplest example, $\msh$ a regular tetrahedron with vertices
\begin{equation}
   A :=
   \left[ \begin{smallmatrix}
      -1 \\ -1\\  -1
   \end{smallmatrix} \right],
   \quad
   B:=
   \left[ \begin{smallmatrix}
      -1 \\  1\\   1
   \end{smallmatrix} \right],
   \quad
   C:=
   \left[ \begin{smallmatrix}
       1 \\ -1\\   1
   \end{smallmatrix} \right],
   \quad
   D:=
   \left[ \begin{smallmatrix}
       1 \\  1\\  -1
   \end{smallmatrix} \right],
\end{equation}
suffices to show that the construction of \cite{Akle:2017:IM}
as stated can not in general generate $G^1$ surfaces.
Let $\bm$ be the point where the curves connecting
the limit points associated with $C$ and $D$ meet
the curves connecting $\bv$, the center of the face $B,C,D$, to 
the center of $A,C,D$ (see \figref{fig:tet}).
We consider $G^1$ continuity along the edge from $\bv$  to $\bm$.

To compute with integers throughout we scale $\msh$ by $2^23^2\cdot5\cdot7$.
Following the algorithm of \cite{Akle:2017:IM} up to the claim
`Our calculation of the control points guarantees $G_1$ continuity',
the mesh points and BB-coefficients can then be computed as integers.
Three rows of BB-coefficients determine the $G^1$ continuity 
constraints \eqref{eq:g1} between the resulting two adjacent bi-3 patches
$\bp$ and $\bq$. 
We focus on on the BB-coefficients of $\bp_{ij}$ for $i=0,1,2,3$ and $j=0,1$
using $\sim$ to indicate proportionality after scaling the 
coefficients to the right of $\sim$ to the smallest integer values:
(after multiplication by 210)
\begin{align}
\begin{matrix}
   \bp_{i1} \\
   \bp_{i0}=\bq_{i,3} \\
   \bq_{i2} 
\end{matrix}
\sim
\begin{matrix}
\left[ \begin{smallmatrix} 7\\10\\7 \end{smallmatrix} \right],
\left[ \begin{smallmatrix} 10\\10\\4 \end{smallmatrix} \right],
\left[ \begin{smallmatrix} 16\\4\\-2 \end{smallmatrix} \right],
\left[ \begin{smallmatrix} 16\\3\\-3 \end{smallmatrix} \right],
\\ 
\left[ \begin{smallmatrix} 8\\8\\8 \end{smallmatrix} \right],
\left[ \begin{smallmatrix} 10\\7\\7 \end{smallmatrix} \right],
\left[ \begin{smallmatrix} 16\\1\\1 \end{smallmatrix} \right],
\left[ \begin{smallmatrix} 16\\0\\0 \end{smallmatrix} \right],
\\ 
\left[ \begin{smallmatrix} 7\\7\\10 \end{smallmatrix} \right],
\left[ \begin{smallmatrix} 10\\4\\10 \end{smallmatrix} \right],
\left[ \begin{smallmatrix} 16\\-2\\4 \end{smallmatrix} \right],
\left[ \begin{smallmatrix} 16\\-3\\3 \end{smallmatrix} \right],
\end{matrix}
\label{eq:g1coef}
\end{align}
 
Then the coefficients of the derivatives across and along the common edge are 
(after multiplication by 630)
\begin{align}
\begin{matrix}
\partial_2 \bp \\
\partial_1 \bp  = \partial_1 \bq\\
\partial_2 \bq \\
\end{matrix}
\sim
\begin{matrix}
&
\left[ \begin{smallmatrix} -1\\2\\-1 \end{smallmatrix} \right], 
\left[ \begin{smallmatrix} 0\\3\\-3 \end{smallmatrix} \right], 
\left[ \begin{smallmatrix} 0\\3\\-3 \end{smallmatrix} \right], 
\left[ \begin{smallmatrix} 0\\3\\-3 \end{smallmatrix} \right], 
\\ 
\quad &
\left[ \begin{smallmatrix} 2\\-1\\-1 \end{smallmatrix} \right], 
\left[ \begin{smallmatrix} 6\\-6\\-6 \end{smallmatrix} \right], 
\left[ \begin{smallmatrix} 0\\-1\\-1 \end{smallmatrix} \right], 
\\ 
&
\left[ \begin{smallmatrix} 1\\1\\-2 \end{smallmatrix} \right], 
\left[ \begin{smallmatrix} 0\\3\\-3 \end{smallmatrix} \right], 
\left[ \begin{smallmatrix} 0\\3\\-3 \end{smallmatrix} \right], 
\left[ \begin{smallmatrix} 0\\3\\-3 \end{smallmatrix} \right]. 
\end{matrix}
\end{align}

We compute the determinant
(scaled by 10716300)
\begin{align}
| \partial_2 \bp, \partial_1 \bp, \partial_2 \bq | 
\sim 
 [ 0,  105,  185,  105,  36,  5,  0,  0,  0].
\end{align}

Taking taking the dot-product of  \eqref{eq:g1}  with 
$\partial_1 \bp(u,0) \times \partial_2 \bq(u,0)$ implies
that \\
$|\partial_2 \bp(u,0), \partial_1 \bp(u,0), \partial_2 \bq(u,0) |=0$.
However in the counterexample the determinant is non-zero and hence
the two patches do not join with $G^1$ continuity.

\section{Alternative bi-3 constructions}
\label{sec:alternative} 
Attempts to generalize bi-cubic splines to irregular layouts have a long 
history and include
\cite{bezier77a},
\cite{beeker86a},
\cite{Catmull-1978-CC},
\cite{Sabin:1968:CCS},
\cite{Sarraga:1987:IGU}, 
\cite{Gregory:1994:FPH},
\cite{Peters:1991:SIM},
\cite{Peters:1994:SAT},
\cite{Wijk:1986:BPA} 
to list just a few.
The $3\times3$ bi-3 patches per quad construction of \cite{Fan:2008:SBS}
achieves the lower bound determined by \thmref{thm:twodblknt} and
has been used to implement the multi-sided caps included in
bi-cubic T-spline constructions.
\cite{Sarov:RPG:2016} focuses on restricted input mesh
to ensure that $G^1$ bi-3 surfaces can be built with fewer pieces.

Use of a guide shape (of higher polynomial degree) appears to be necessary
to construct bi-3 surfaces with a good distribution of highlight lines,
as required for car styling and many other outer surfaces. 
For example, the guided approach improves the shape 
of bi-3 singularly parameterized surfaces \cite{Karciauskas:CCB:2016}.
The paper ``Can bi-cubic surfaces be class {A}?''
\cite{Karciauskas:2015:CBSconf} emphasizes the distinction
between exact $G^1$ continuity and acceptable shape
in terms of curvature distribution and highlight lines.
This distinction, accompanied by mathematical estimates of the 
jump in normals, could also be useful in the context of 
\cite{Akle:2017:IM}. Since proving surface `fairness' is typically not possible,
it is recommended to test new surface construction algorithms
on the obstacle course \cite{obstaclecourse} of local input meshes.


\section{Conclusion}
The approach of \cite{Akle:2017:IM} rounds shapes 
but cannot guarantee $G^1$ continuity.
A number of alternative finite bi-3 surface constructions exist
in the literature. Depending on the valence they
require more or fewer pieces than \cite{Akle:2017:IM}.
There are many constructions
using few patches but of higher degree than bi-3.


\bibliographystyle{alpha}
\bibliography{p,/cise/homes/jorg/public_html/peters-jorg}
\end{document}